\documentclass[a4paper,english]{article}
\usepackage{ifthen}
\usepackage[english]{babel}
\usepackage{bibgerm}
\usepackage[utf8]{inputenc}
\usepackage{graphicx} 
\usepackage[hyphens]{url}
\usepackage{hyperref}
\usepackage{caption}
\usepackage{amsmath}
\usepackage{amsthm}
\usepackage{amsfonts}
\usepackage{amssymb}
\usepackage{mathtools}
\usepackage{enumerate}
\usepackage{xcolor}
\usepackage{xspace}
\usepackage{times}
\usepackage[a4paper, total={6in, 9in}]{geometry}
\usepackage{natbib}

\allowdisplaybreaks

\usepackage{zibtitlepage}
\usepackage{authblk}

\usepackage{algorithm}
\usepackage{algpseudocode}

\begin{document}

\title{GPU Parallel Pareto Local Search for Multi-objective QUBO}

\author[1]{Stephanie Riedmüller\footnote{corresponding author, riedmueller@zib.de, \ZTPOrcid{0009-0006-4508-4262}}
}
\author[1,2]{Thorsten Koch
\footnote{\ZTPOrcid{0000-0002-1967-0077}}
}
\affil[1]{Zuse Institute Berlin, Berlin, Germany}
\affil[2]{Technische Universität Berlin, Berlin, Germany}

\newcommand{\titlepageauthors}{Stephanie Riedmüller$^{1}$, Thorsten Koch$^{1,2}$}
\newcommand{\titlepageaffiliations}{
$^{1}$Zuse Institute Berlin, Berlin, Germany\\
$^{2}$Technische Universtität Berlin, Berlin, Germany
}

\newcommand{\acknowledgments}{
The work for this article has been conducted in the Research Campus MODAL funded by the Federal Ministry of Research, Technology and Space (BMFTR) (fund numbers 05M14ZAM, 05M20ZBM, 05M2025).
}

\maketitle

\begin{abstract}
We revisit Pareto Local Search (PLS) for multi-objective quadratic unconstrained binary optimization (QUBO) problems. Multi-objective QUBO formulations have gained increasing attention as a universal modeling method and natural input for quantum optimization. While exact methods exist, their unconstrained binary structure makes them well-suited for heuristics such as PLS. At the same time, GPUs offer a powerful classical alternative to quantum-inspired optimization. We develop an efficient GPU-based parallel implementation of PLS for multi-objective QUBO problems. Our approach combines coarse-grained and fine-grained parallelization. On the coarse-grained level, an initial set of supported solutions is generated using weighted-sum optimization, from which multiple PLS processes are started in parallel. At the fine-grained level, neighborhood exploration is parallelized on the GPU by simultaneously evaluating binary-flip candidate solutions across all objectives. Non-dominated solutions are maintained using an ND-Tree data structure to support efficient Pareto updates. The implementation combines several known PLS techniques within a unified CPU/GPU hybrid architecture. We describe the algorithmic components, the GPU-specific design decisions, and report results from a computational study evaluating the effectiveness of the proposed approach.
\end{abstract}

\section{Introduction}
\label{sec:intro}
The effectiveness of an optimization algorithm strongly depends on the problem class. While some combinatorial optimization problems can be solved efficiently with exact methods, unconstrained problems are particularly well suited to sampling-based heuristics.
One important class of unconstrained problems is Quadratic Unconstrained Binary Optimization (QUBO), which has gained increasing attention due to its applicability to quantum computing. It seems to be the same property that makes them interesting for quantum computing as well as classical heuristics, as both rely on sampling. This paper investigates how such classical heuristics can be further accelerated using specialized hardware, namely GPUs. To that end, we revisit Pareto Local Search (PLS) and combine several established techniques into a new GPU parallel PLS implementation.

Given $n\in \mathbb{N}$, $p\geq 2$, and objective functions $f_1, \ldots, f_p: \{0,1\}^n \rightarrow \mathbb{R}$ with $f = (f_1, \ldots, f_p)$ and $f_i = x^T Q_i x$ for $i=1, \ldots p$ and $Q_i \in \mathbb{R}^{n \times n}$, we consider the multi-objective Quadratic Unconstrained Binary Optimization problem
\begin{align}
    \min_{x \in \{0,1\}^n} \quad &(f_1(x), \ldots, f_p(x)). \tag{mQUBO} \label{eq:mqubo} 
\end{align}
This problem is also sometimes referred to as multi-objective Unconstrained Binary Quadratic Programming (mUBQP).
We call $x \in \{0,1\}^n$ a \emph{solution} to  \ref{eq:mqubo} in the \emph{decision space} $X = \{0,1\}^n$ and $y = f(x)$ a \emph{point} in the \emph{objective space} $f(X)$.    
A point $\hat{y} \in f(X)$ \emph{dominates} a point $y \in f(X)$, if $\hat{y}_i \leq y_i$ for all $i \in \{1, \ldots, p\}$ with $\hat{y}_i < y_i$ for at least one $i$. 
A solution $\hat{x}\in X$ to  \ref{eq:mqubo} is called \emph{efficient} or \emph{Pareto-optimal}, if the corresponding point in the objective space is not dominated by any other point. In this case, $f(\hat{x})$ is called \emph{non-dominated point}.
The set of efficient solutions is called the \emph{efficient set} and the corresponding set of non-dominated points is called the \emph{non-dominated set} or \emph{Pareto front}. 
For a general introduction to multi-objective optimization, we refer to \citet{Ehrgott2005}.
When comparing the results of multi-objective optimization algorithms, the hypervolume indicator is one of the most widely used quality measures \citep{Riquelme2015}. 

PLS is the natural extension of single-objective Local Search (LS). The goal of PLS is to find (a subset of) the Pareto front by iteratively exploring the neighborhood of already found solutions in the decision space. Therefore, classical PLS involves the following steps: initialization, solution selection, and neighborhood exploration. A generic PLS is given in Algorithm~\ref{alg:pls}.

\begin{algorithm}
\caption{Pareto Local Search}
\label{alg:pls}
\begin{algorithmic}[1]
\Require Initial (solution, point) pair $(s_0, f(s_0))$
\Ensure Pareto front $PF$
\State Initialize $PF \gets \{(s_0, f(s_0))\}$ 
\While{stopping criteria not met}
    \State Select an unexplored solution $s \in PF$ and mark it as explored
    \State Explore the neighborhood $N(s)$ by computing the objective values $f(s')$ for all $s' \in N(s)$
    \State Update $PF$ by keeping only (efficient solution, non-dominated point) pairs
\EndWhile
\State \Return $PF$
\end{algorithmic}
\end{algorithm}

\section{Literature Review}
\label{sec:lit}
For the single-objective case, LS has been intensively researched and there exists a parallel GPU version \citep{Luong2013}. 
In the multi-objective setting, for classical sequential PLS, the latest research leads to an acceleration in the neighborhood exploration step \citep{Kang2024}, where different strategies are introduced and evaluated. Best improvements could be made by the so-called Don't look bits, where bits are marked as uninteresting, if they did not lead to improvements before. 
PLS has been parallelized on CPU based on decomposition (PPLS/D) \citep{Shi2018}. The search space is decomposed into cones and each cone is searched in parallel on the CPU. Instead of Pareto dominance, the search starts as single-objective LS with a weighted Tchebycheff objective to direct the search along the cone. Later, the search in each cone follows with classical PLS based on Pareto dominance. Sharing solutions between the parallel search cones and redefining the decomposition of uneven search spaces at fixed time steps provides only marginal improvements to PPLS/D \citep{Shi2024}.
Comparable to the single-objective LS, there also exists a Tabu Search version for PLS, the so-called Pareto Tabu Search (PTS).
For continuous optimization, a parallel PTS has been published, which performs the m objective function evaluations in parallel on the CPU and also introduces several memory structures (short-term, medium-term, long-term) for Tabu solutions, intensification solutions and the Pareto front, respectively \citep{Jaeggi2008}. A corresponding parallel version on GPU is also described for continuous problems by \citet{Tsotskas2014}.
Apart from PLS variants, similar techniques are used in many variants of parallel evolutionary algorithms both on CPU and GPU, mostly based on the popular algorithms NSGA-II and MOEA/D. One of the latest approaches is given by \citet{Sun2024}, which introduces a hybrid CPU/GPU scheme for NSGA-II combined with MOEA/D. 
Independent of the individual algorithm, storing and updating the current Pareto front (dynamic non-dominance problem) is usually a huge bottleneck. Next to simple lists for few objectives, a successful approach for at least 3 objectives is to use ND-Trees \citep{Jaszkiewicz2018}.

\section{Parallel PLS and GPU Implementation}
\label{sec:3}

For the following algorithmic choices, we combine several established techniques for single-objective LS and sequential and parallel PLS (PPLS).
The main structure of the algorithm is based on PPLS/D \citep{Shi2018}.
Since we aim at using the hardware to its fullest potential, the goal is to have minimal data transfer between CPU and GPU and an algorithmic design in line with the SIMD principle.

\subsubsection*{Initialization / parallelization on algorithmic level}

We start PLS for a set of different start solutions in parallel. In single-objective LS, starting from several solutions is usually done to get a robust solution, however, in PLS, this is done to access different areas of the Pareto front. This first set of solutions can either be generated randomly or by computing some supported solutions via a weighted sum approach. We follow the approach given in PPLS/D to decompose the objective space into subregions and use single-objective LS with a weighted sum scalarization to find spread-out start solutions \citep{Shi2018}. Each subregion is treated by a separate worker, a CPU thread or a GPU thread block, respectively.

\subsubsection*{Solution selection}
In contrast to single-objective LS, PLS has to select the current solution out of several non-dominated solutions to explore its neighborhood.  
On the CPU, each subregion selects it by Tchebycheff scalarization as suggested by \citet{Shi2018}. On the GPU, each weight vector first runs a complete Tchebycheff LS on-device, and the resulting local optima are then refined by PLS on the shared archive.

\subsubsection*{Neighborhood exploration / parallelization on iteration level}

The algorithm performs a parallel evaluation of the 1 hamming distance neighborhood. That is, for $n$ 1-flip neighbors, $m$ objective evaluations are executed on the GPU in parallel.
We include the binary representation mapping suggested for the single-objective LS by \citet{Luong2013} and compute the binary flips directly on the GPU to minimize the data transfer between CPU and GPU. We consider all possible binary flips as it is GPU-friendly to consider tasks with minimal if-statements. 
The computation of the objective values follows the incremental QUBO update suggested by \citet{Nakano2023}.
The new solutions are then inserted into the Pareto front archive, which is implemented as an ND-Tree \citep{Jaszkiewicz2018}.

\section{Computational Results}
\label{sec:4}

The mQUBO instances were generated by the mUBQP generator in the mocobench available at \url{https://sourceforge.net/projects/mocobench/}. The instances vary in the number of objectives $m=3,4,5$ and the number of variables $n=100,500,1000$.
For the given instances and a time limit of 1000s, we compare the anytime behavior regarding the hypervolume indicator for three algorithms:
\begin{itemize}
    \item The original PPLS/D algorithm \citep{Shi2018, Shi2024} as the state of the art parallel PLS method,
    \item our new PPLS-cpu implementation exclusively using the CPU backend, and
    \item our new PPLS-gpu implementation using the GPU backend.
\end{itemize}
The two new proposed methods PPLS-cpu and PPLS-gpu are implemented in Rust.\footnote{The implementation and the instances are available at \url{https://git.zib.de/sriedmueller/ppls-cpu-gpu}.}
All computations have been carried out on 
an AMD Ryzen Threadripper PRO 7965WX using 32 cores and a  
NVIDIA GeForce RTX 4090 GPU.

\begin{figure}
    \centering
    \includegraphics[width=\linewidth, trim={0 0 0 1cm}, clip]{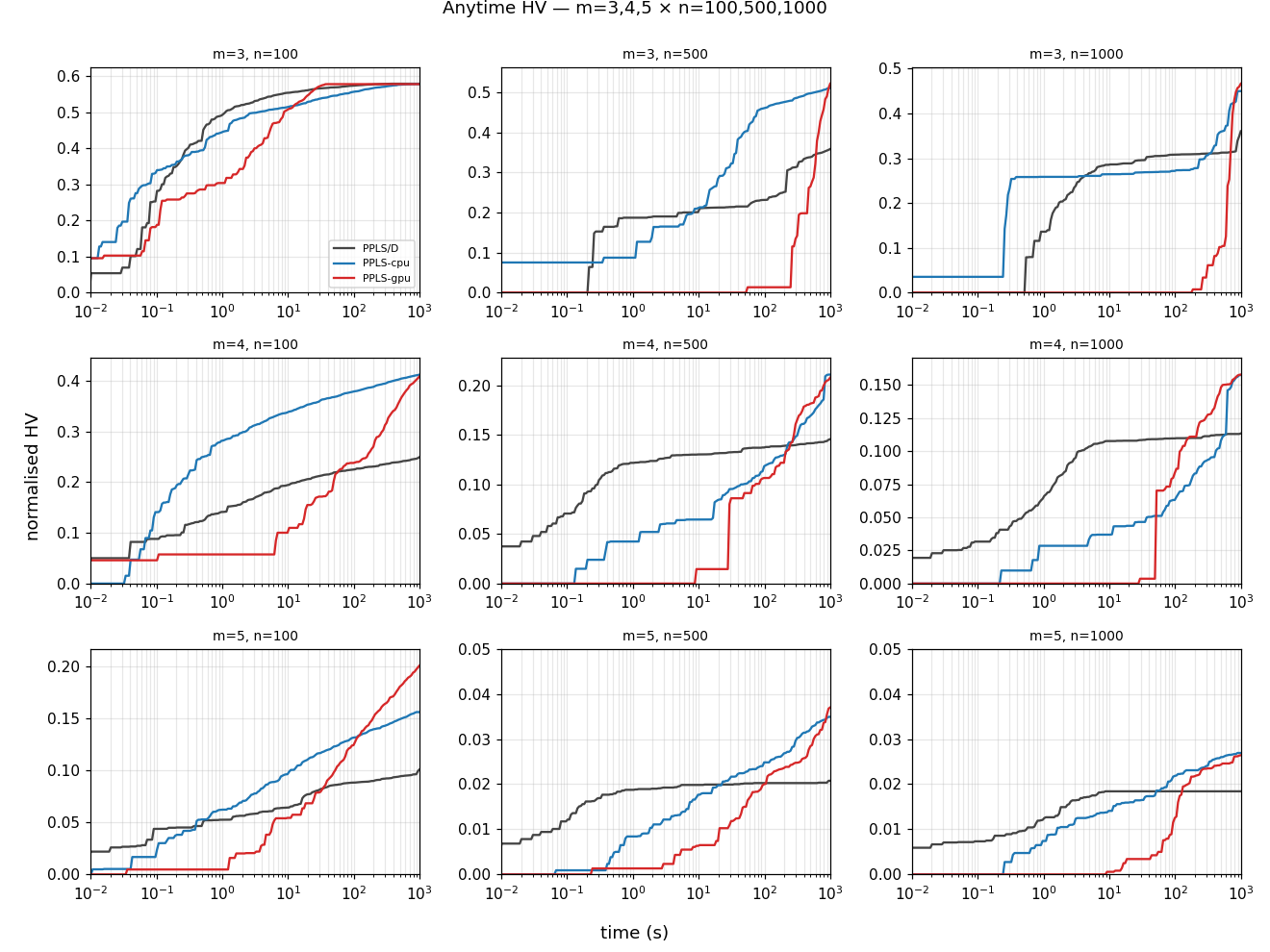}
    \caption{Anytime behavior of the PPLS-cpu, PPLS-gpu and PPLS/D for different instance sizes and number of objectives.}
    \label{fig:results}
\end{figure}

Figure~\ref{fig:results} shows the so-called anytime behavior of the three algorithms on the considered instances, i.e., the normalized hypervolume over time. The hypervolume has been measured via Monte-Carlo estimation \citep{Bringmann_2010}. 
When comparing the anytime-behavior of PPLS/D and PPLS-cpu, the first part of the solving time is dominated by PPLS/D for 4 and 5 objectives increasingly with instance size. The PPLS/D hypervolume curve then flattens and is dominated by PPLS-cpu, which shows a steeper ascent. For 3 objectives, this pattern is not that prominently visible.
PPLS-gpu shows the steepest ascent in hypervolume.
Both PPLS-cpu and PPLS-gpu start later than PPLS/D to increase the hypervolume for two reasons: First, while PPLS/D with default parameters starts with a random solution set, both PPLS-cpu and PPLS-gpu generate the start solutions with a scalarization approach. 
Second, PPLS-gpu completes the Tchebycheff LS of all weight vectors on the GPU before the first solutions enter the archive and, therefore, starts even later than PPLS-cpu.
In the time limit, PPLS-cpu and PPLS-gpu result in a hypervolume up to two times larger than PPLS/D. Note the percentage of hypervolume reached by all three algorithms decreases strongly with number of objectives and instance size, as expected.
Due to the later start of PPLS-gpu and its steeper increase in hypervolume, an even larger performance advantage may be expected for longer runtimes.

\section{Conclusion}
\label{sec:5}

We presented two parallel variants for PLS based on PPLS/D that incorporate ND-Trees and the incremental QUBO update: PPLS-cpu executes the parallelization exclusively on the CPUs, while PPLS-gpu uses the GPU. In a computational study, we showed that both PPLS-cpu and PPLS-gpu beat PPLS/D in the overall hypervolume reached in the time limit. The steepest increase in hypervolume was observed for PPLS-gpu. The performance advantage of PPLS-gpu over PPLS-cpu is rather small in comparison to other examples of GPU-acceleration, which can be explained by the Pareto front archive being the main bottleneck, which is executed on the CPU in both algorithms.

\section{Acknowledgments}
The work for this article has been conducted in the Research Campus MODAL funded by the Federal Ministry of Research, Technology and Space (BMFTR) (fund numbers 05M14ZAM, 05M20ZBM, 05M2025). 

The implementation of PPLS-cpu and PPLS-gpu was developed with the assistance of the AI coding tool Claude Code (Anthropic). All algorithmic design decisions, the experimental setup, and the evaluation of results were made and verified by the authors.

\bibliographystyle{apalike}
\bibliography{references}

@article{Kang2024,
title = {New techniques to improve neighborhood exploration in pareto local search},
journal = {Expert Systems with Applications},
volume = {254},
pages = {124296},
year = {2024},
issn = {0957-4174},
doi = {https://doi.org/10.1016/j.eswa.2024.124296},
url = {https://www.sciencedirect.com/science/article/pii/S095741742401162X},
author = {Yuhao Kang and Jialong Shi and Jianyong Sun and Qingfu Zhang and Ye Fan}
}

@article{Sun2024,
  author    = {Zhuoran Sun and Ying Ying Liu and Parimala Thulasiraman},
  title     = {Cooperative, collaborative, coevolutionary multi-objective optimization on CPU-GPU multi-core},
  journal   = {The Journal of Supercomputing},
  year      = {2024},
  volume    = {81},
  number    = {1},
  pages     = {59},
  date      = {2024-10-22},
  issn      = {1573-0484},
  doi       = {10.1007/s11227-024-06525-8},
  url       = {https://doi.org/10.1007/s11227-024-06525-8}
}

@inproceedings{Nakano2023,
   title={Diverse Adaptive Bulk Search: a Framework for Solving QUBO Problems on Multiple GPUs},
   url={http://dx.doi.org/10.1109/IPDPSW59300.2023.00060},
   DOI={10.1109/ipdpsw59300.2023.00060},
   booktitle={2023 IEEE International Parallel and Distributed Processing Symposium Workshops (IPDPSW)},
   publisher={IEEE},
   author={Nakano, Koji and Takafuji, Daisuke and Ito, Yasuaki and Yazane, Takashi and Yano, Junko and Ozaki, Shiro and Katsuki, Ryota and Mori, Rie},
   year={2023},
   month=May, pages={314–325} }

@article{Tsotskas2014,
title = {The Design and Implementation of a GPU-enabled Multi-objective Tabu-search Intended for Real World and High-dimensional Applications},
journal = {Procedia Computer Science},
volume = {29},
pages = {2152-2161},
year = {2014},
note = {2014 International Conference on Computational Science},
issn = {1877-0509},
doi = {https://doi.org/10.1016/j.procs.2014.05.200},
url = {https://www.sciencedirect.com/science/article/pii/S1877050914003779},
author = {Christos Tsotskas and Timoleon Kipouros and Anthony Mark Savill}}

@article{Shi2024,
  author    = {Shi, J. and Sun, J. and Zhang, Q. and Zhang, H. and Fan, Y.},
  title     = {Improving Pareto Local Search Using Cooperative Parallelism Strategies for Multiobjective Combinatorial Optimization},
  journal   = {IEEE Transactions on Cybernetics},
  year      = {2024},
  volume    = {54},
  number    = {4},
  pages     = {2369--2382},
  doi       = {10.1109/TCYB.2022.3226744},
  pmid      = {37015461},
  note      = {Epub 2024 Mar 18}
}

@article{Shi2018,
  author={Shi, Jialong and Zhang, Qingfu and Sun, Jianyong},
  journal={IEEE Transactions on Cybernetics}, 
  title={PPLS/D: Parallel Pareto Local Search Based on Decomposition}, 
  year={2020},
  volume={50},
  number={3},
  pages={1060-1071},
  doi={10.1109/TCYB.2018.2880256}}

@article{Jaszkiewicz2018,
author = {Jaszkiewicz, Andrzej and Lust, Thibaut},
title = {ND-Tree-Based Update: A Fast Algorithm for the Dynamic Nondominance Problem},
year = {2018},
issue_date = {Oct. 2018},
publisher = {IEEE Press},
volume = {22},
number = {5},
issn = {1089-778X},
url = {https://doi.org/10.1109/TEVC.2018.2799684},
doi = {10.1109/TEVC.2018.2799684},
journal = {Trans. Evol. Comp},
month = oct,
pages = {778–791},
numpages = {14}
}

@ARTICLE{Luong2013,
  author={Luong, Thé Van and Melab, Nouredine and Talbi, El-Ghazali},
  journal={IEEE Transactions on Computers}, 
  title={GPU Computing for Parallel Local Search Metaheuristic Algorithms}, 
  year={2013},
  volume={62},
  number={1},
  pages={173-185},
  doi={10.1109/TC.2011.206}}

@article{Jaeggi2008,
title = {The development of a multi-objective Tabu Search algorithm for continuous optimisation problems},
journal = {European Journal of Operational Research},
volume = {185},
number = {3},
pages = {1192-1212},
year = {2008},
issn = {0377-2217},
doi = {https://doi.org/10.1016/j.ejor.2006.06.048},
url = {https://www.sciencedirect.com/science/article/pii/S0377221706006357},
author = {D.M. Jaeggi and G.T. Parks and T. Kipouros and P.J. Clarkson}
}

@INPROCEEDINGS{Riquelme2015,
  author={Riquelme, Nery and Von Lücken, Christian and Baran, Benjamin},
  booktitle={2015 Latin American Computing Conference (CLEI)}, 
  title={Performance metrics in multi-objective optimization}, 
  year={2015},
  volume={},
  number={},
  pages={1-11},
  doi={10.1109/CLEI.2015.7360024}}

@book{Ehrgott2005,
    author = {Matthias Ehrgott},
    title = {Multicriteria Optimization},
    publisher = {Springer Berlin, Heidelberg},
    year = {2005},
    DOI = {https://doi.org/10.1007/3-540-27659-9}
}

@article{Bringmann_2010,
   title={Approximating the volume of unions and intersections of high-dimensional geometric objects},
   volume={43},
   ISSN={0925-7721},
   url={http://dx.doi.org/10.1016/j.comgeo.2010.03.004},
   DOI={10.1016/j.comgeo.2010.03.004},
   number={6–7},
   journal={Computational Geometry},
   publisher={Elsevier BV},
   author={Bringmann, Karl and Friedrich, Tobias},
   year={2010},
   month=aug, pages={601–610}
}

\end{document}